\documentclass[12pt]{article}
\usepackage{amssymb,amsmath,amscd, amsthm,stmaryrd,verbatim, mathrsfs}

\begin{document}
\baselineskip=14pt

\newcommand{\la}{\langle}
\newcommand{\ra}{\rangle}
\newcommand{\psp}{\vspace{0.4cm}}
\newcommand{\pse}{\vspace{0.2cm}}
\newcommand{\ptl}{\partial}
\newcommand{\dlt}{\delta}
\newcommand{\sgm}{\sigma}
\newcommand{\al}{\alpha}
\newcommand{\be}{\beta}
\newcommand{\G}{\Gamma}
\newcommand{\gm}{\gamma}
\newcommand{\vs}{\varsigma}
\newcommand{\Lmd}{\Lambda}
\newcommand{\lmd}{\lambda}
\newcommand{\td}{\tilde}
\newcommand{\vf}{\varphi}
\newcommand{\yt}{Y^{\nu}}
\newcommand{\wt}{\mbox{wt}\:}
\newcommand{\rd}{\mbox{Res}}
\newcommand{\ad}{\mbox{ad}}
\newcommand{\stl}{\stackrel}
\newcommand{\ol}{\overline}
\newcommand{\ul}{\underline}
\newcommand{\es}{\epsilon}
\newcommand{\dmd}{\diamond}
\newcommand{\clt}{\clubsuit}
\newcommand{\vt}{\vartheta}
\newcommand{\ves}{\varepsilon}
\newcommand{\dg}{\dagger}
\newcommand{\tr}{\mbox{Tr}}
\newcommand{\ga}{{\cal G}({\cal A})}
\newcommand{\hga}{\hat{\cal G}({\cal A})}
\newcommand{\Edo}{\mbox{End}\:}
\newcommand{\for}{\mbox{for}}
\newcommand{\kn}{\mbox{ker}}
\newcommand{\Dlt}{\Delta}
\newcommand{\rad}{\mbox{Rad}}
\newcommand{\rta}{\rightarrow}
\newcommand{\mbb}{\mathbb}
\newcommand{\lra}{\Longrightarrow}
\newcommand{\X}{{\cal X}}
\newcommand{\Y}{{\cal Y}}
\newcommand{\Z}{{\cal Z}}
\newcommand{\U}{{\cal U}}
\newcommand{\V}{{\cal V}}
\newcommand{\W}{{\cal W}}
\newcommand{\sta}{\theta}
\setlength{\unitlength}{3pt}
\newcommand{\msr}{\mathscr}

\begin{center}{\Large \bf Projective Oscillator Representations \\ of $sl(n+1)$ and $sp(2m+2)$} 
\end{center}
\vspace{0.2cm}

\begin{center}{\large Xiaoping Xu
}\footnote{Research supported
 by NSFC Grants 11171324 and  11321101.}\end{center}
\begin{center}{Hua Loo-Keng Key Mathematical Laboratory\\
Institute of Mathematics, Academy of Mathematics \& System
Sciences\\ Chinese Academy of Sciences, Beijing 100190, P.R. China\\
Email: xiaoping@math.ac.cn. Tel: 86-10-62651308. Fax: 86-10-62553022
}\end{center}

\begin {abstract}
\quad

The $n$-dimensional projective group gives rise to a one-parameter
family of inhomogeneous first-order differential operator
representations of $sl(n+1)$. By partially swapping differential
operators and multiplication operators, we obtain more general
differential operator  representations of $sl(n+1)$. Letting these
differential operators act on the corresponding polynomial algebra
and the space of exponential-polynomial functions, we construct new
multi-parameter families of explicit infinite-dimensional
irreducible representations for $s(n+1)$ and $sp(2m+2)$ when
$n=2m+1$. Our results can be viewed as  extensions of Howe's
oscillator construction of infinite-dimensional multiplicity-free
irreducible representations for $sl(n)$.\vspace{0.3cm}

\noindent{\it Keywords}:\hspace{0.3cm} special linear Lie algebra;
symplectic Lie algebra; oscillator
 representation; irreducible module; polynomial algebra; exponential-polynomial function.

\end{abstract}

\section {Introduction}

\quad$\;$ A module of a finite-dimensional simple Lie algebra is
called a {\it weight module} if it is a direct sum of its weight
subspaces. A module of a finite-dimensional simple Lie algebra is
called {\it cuspidal} if it is not induced from its proper parabolic
subalgebras. Infinite-dimensional irreducible weight modules of
finite-dimensional simple Lie algebras with finite-dimensional
weight subspaces have been intensively studied by the authors in
\cite{BBL, BFL, BHL, BL1, BL2, Fs, Fv, M}. In particular, Fernando
\cite{Fs} proved that such modules must be cuspidal or parabolically
induced. Moreover, such cuspidal modules exist only for special
linear Lie algebras and symplectic Lie algebras. A similar result
was independently obtained by Futorny \cite{Fv}. Mathieu \cite{M}
proved that such cuspidal  modules
 are irreducible components in the tensor
modules of their multiplicity-free modules with finite-dimensional
modules. Although the structures of irreducible weight modules of
finite-dimensional simple Lie algebras with finite-dimensional
weight subspaces were essentially determined by Fernando's result in
\cite{Fs} and Methieu's result in \cite{M}, explicit structures of
such modules are not that known. It is important to find explicit
natural realizations of them.

Let $\mbb{F}$ be a field with characteristic $0$ (say,
$\mbb{Q},\mbb{R},\mbb{C}$) and let $n\geq 2$ be an integer. A
projective transformation on $\mathbb{F}^n$ is given by
$$u\mapsto \frac{Au+\vec b}{\vec c\:^t u+d}\qquad\mbox{for}\;\;u\in
\mathbb{F}^n,\eqno(1.1)$$ where all the vectors in $\mathbb{F}^n$
are in column form and
$$\left(\begin{array}{cc}A&\vec b\\ \vec c\:^t&
d\end{array}\right)\in GL(n).\eqno(1.2)$$ It is well-known that a
transformation of mapping straight lines to lines must be a
projective transformation. The above transformations give rise to an
inhomogeneous representation of the Lie algebra $sl(n+1,\mbb{F})$ on
the polynomial functions of the projective space. Using Shen's mixed
product for Witt algebras in \cite{S} and the above representation,
Zhao and the author \cite{ZX}
 constructed a new functor from
$gl(n,\mbb{F})$-{\bf Mod} to $sl(n+1)$-{\bf Mod} and found a
condition for the functor to map a finite-dimensional irreducible
$gl(n,\mbb{F})$-module to an infinite-dimensional irreducible
$sl(n+1,\mbb{F})$-module. Our general frame also gave a direct
polynomial extension from irreducible $gl(n,\mbb{F})$-modules to
irreducible $sl(n+1,\mbb{F})$-modules.

The work \cite{ZX} lead to a one-parameter family of inhomogeneous
first-order differential operator (oscillator) representations of
$sl(n+1,\mbb{F})$. By partially swapping differential operators and
multiplication operators, we obtain more general differential
operator (oscillator) representations of $sl(n+1,\mbb{F})$. In this
paper, we construct new multi-parameter families of explicit
infinite-dimensional irreducible representations for
$s(n+1,\mbb{F})$ and $sp(2m+2)$ when $n=2m+1$ by letting these
differential operators act on the corresponding polynomial algebra
and the space of exponential-polynomial functions. Some of the
corresponding modules are explicit infinite-dimensional irreducible
weight modules with finite-dimensional weight subspaces. Our results
can be viewed as extensions of Howe's oscillator construction of
infinite-dimensional multiplicity-free irreducible representations
for $sl(n,\mbb{F})$ (cf. \cite{H}). Indeed, Howe's result plays an
important role in proving the irreducibility of the representations
for $sl(n+1,\mbb{F})$. The results on symplectic Lie algebras in
this paper can be used to study the irreducible representations of
the other simple Lie algebra via Howe's theta correspondence
technique.

 Let $E_{r,s}$ be the $(n+1)\times(n+1)$
matrix with 1 as its $(r,s)$-entry and 0 as the others. The special
linear algebra
$$sl(n+1,\mbb{F})=\sum_{1\leq i<j\leq
n+1}(\mbb{F}E_{i,j}+\mbb{F}E_{j,i})+\sum_{r=1}^n\mbb{F}(E_{r,r}-E_{r+1,r+1}).\eqno(1.3)$$
 For any
two integers $p\leq q$, we denote $\ol{p,q}=\{p,p+1,\cdots,q\}$. Set
$D=\sum_{s=1}^nx_s\ptl_{x_s}$. According to Zhao and the author's
work \cite{ZX},  we have the following one-parameter generalization
$\pi_c$ of the projective representation  of $sl(n+1,\mbb{F})$:
$$\pi_c(E_{i,j})=x_i\ptl_{x_j},\;\;\pi_c(E_{i,n+1})=x_i(D+c),\;\;\pi_c(E_{n+1,i})=-\ptl_{x_i},\eqno(1.4)$$
$$\pi_c(E_{i,i}-E_{j,j})=x_i\ptl_{x_i}-x_j\ptl_{x_j},\;\;\pi_c(E_{n,n}-E_{n+1,n+1})=D+c+x_n\ptl_{x_n}\eqno(1.5)$$
 for $i,j\in\ol{1,n}$ with $i\neq j$, where $c\in\mbb{F}$.

Let $S$ be a subset of $\ol{1,n}$. Note the symmetry:
$$[\ptl_{x_r},x_r]=1=[-x_r,\ptl_{x_r}].\eqno(1.6)$$
Changing operators $\ptl_{x_r}\mapsto -x_r$ and $x_r\mapsto
\ptl_{x_r}$ for $r\in S$ in (1.4) and (1.5), we get another
differential-operator representation $\pi_{c,S}$ of
$sl(n+1,\mbb{F})$. We treat $\pi_{c,\emptyset}=\pi_c$ and call
$\pi_{c,S}$ {\it  projective oscillator representations} in terms of
physics terminology.
 For
$\vec a=(a_1,a_2,...,a_n)^t\in\mbb{F}^n$, we denote $\vec a\cdot\vec
x=\sum_{i=1}^na_ix_i$. Let ${\msr A}=\mbb{F}[x_1,x_2,...,x_n]$ be
the algebra of polynomials in $x_1,x_2,...,x_n$. Moreover, we set
$${\msr A}_{\vec a}=\{fe^{\vec a\cdot\vec
x}\mid f\in{\msr A}\}.\eqno(1.7)$$ Denote by $\pi_{c,S}^{\vec a}$
the representation $\pi_{c,S}$ of $sl(n+1,\mbb{F})$ on ${\msr
A}_{\vec a}$ and by $\mbb{N}$ the set of nonnegative integers. In
\cite{ZX}, Zhao and the author proved that the representation
$\pi_{c,\emptyset}^{\vec 0}$ of $sl(n+1,\mbb{F})$ is irreducible  if
and only if $c\not\in -\mbb{N}$. Moreover, ${\msr A}$ has a
composite series of length 2 when $c\in -\mbb{N}$. In this paper, we
prove:\psp

{\bf Theorem 1}. {\it Let $S$ be a proper subset of $\ol{1,n}$. The
representation $\pi_{c,S}^{\vec 0}$ is irreducible for any $c\in\mbb
F\setminus \mbb Z$, and the underlying module ${\msr A}$ is an
infinite-dimensional weight $sl(n+1,\mbb{F})$-module with
finite-dimensional weight subspaces. If $a_i\neq 0$ for some $i\in
\ol{1,n}\setminus S$ or $|S|>1$ and $\vec a\neq 0$, then the
representation $\pi_{c,S}^{\vec a}$ of $sl(n+1,\mbb{F})$ is always
irreducible for any $c\in\mbb F$.} \psp

Suppose that $n=2m+1>1$ is an odd integer and the subset $S$
satisfies:
$$m+1\not\in S\;\;\mbox{and for}\;i\in\ol{1,m},\;\mbox{at most one
of}\;i\;\mbox{and}\;i+m+1\;\mbox{in}\; S.\eqno(1.8)$$ Our second
main theorem in this paper is as follows.\psp

{\bf Theorem 2}. {\it If $c\not\in -\mbb{N}$, the restricted
representation $\pi_{c,\emptyset}^{\vec 0}$ of $sp(2m+2,\mbb{F})$ is
irreducible. When $c\in -\mbb{N}$, the $sp(2m+2,\mbb{F})$-module
${\msr A}$ has a composite series of length 2 with respect to the
restricted representation $\pi_{c,\emptyset}^{\vec 0}$.

The restricted representation $\pi_{c,S}^{\vec 0}$ of
$sp(2m+2,\mbb{F})$ with $S\neq\emptyset$ is irreducible for any
$c\in\mbb{F}\setminus \mbb Z$. Suppose that $\vec a\neq\vec
0,\;a_{m+1}=0, \;a_{i_0}\neq 0$ for some $m+1+i_0\in
S\bigcap\ol{m+2,2m+1}$ if $S\bigcap\ol{m+2,2m+1}\neq\emptyset$, and
$a_{m+1+j_0}\neq 0$ for some $j_0\in S\bigcap\ol{1,m+1}$ if
$S\bigcap\ol{1,m+1}\neq\emptyset$, then the restricted
representation $\pi_{c,S}^{\vec a}$ of $sp(2m+2,\mbb{F})$ is
irreducible for any $c\in\mbb{F}$.

With respect to the restricted representation $\pi_{c,S}^{\vec 0}$,
${\msr A}$ is an infinite-dimensional weight
$sp(2m+2,\mbb{F})$-module with finite-dimensional weight subspaces.
}\psp

In Section 2, we prove Theorem 1. The proof of Theorem 2 is given in
Section 3.

\section{Proof of Theorem 1}

In this section, we will prove Theorem 1 case by case.\psp

{\it Case 1. The representation $\pi_{c,S}^{\vec 0}$ with $S\neq
\emptyset,\ol{1,n}$.} \psp

Without loss of generality, we assume $S=\ol{1,n_1}$ for some
$n_1\in\ol{1,n_1}$ and $n_1<n$.  Set
$$\td D=\sum_{r=n_1+1}^nx_r\ptl_{x_r}-\sum_{i=1}^{n_1}x_i\ptl_{x_i}.\eqno(2.1)$$
Then the representation $\pi_{c,S}^{\vec 0}$ of $sl(n+1,\mbb{F})$ is
the representation $\pi_{c,S}$ on ${\msr A}$ with
$$\pi_{c,S}(E_{i,j})=\left\{\begin{array}{ll}-x_j\ptl_{x_i}-\delta_{i,j}&\mbox{if}\;
i,j\in\ol{1,n_1};\\ \ptl_{x_i}\ptl_{x_j}&\mbox{if}\;i\in\ol{1,n_1},\;j\in\ol{n_1+1,n};\\
-x_ix_j &\mbox{if}\;i\in\ol{n_1+1,n},\;j\in\ol{1,n_1};\\
x_i\partial_{x_j}&\mbox{if}\;i,j\in\ol{n_1+1,n},
\end{array}\right.\eqno(2.2)$$
$$\pi_{c,S}(E_{i,n+1})=\left\{\begin{array}{ll}(\td D+c-n_1-1)\ptl_{x_i}&\mbox{if}\;\;i\leq
n_1,\\ \\
x_i(\td D+c-n_1)&\mbox{if}\;\;i> n_1,\end{array}\right.\eqno(2.3)$$
$$\pi_{c,S}(E_{n+1,i})=\left\{\begin{array}{ll}x_i&\mbox{if}\;\;i\leq
n_1,\\-\ptl_{x_i}&\mbox{if}\;\;i> n_1,\end{array}\right.\eqno(2.4)$$
$$\pi_{c,S}(E_{n,n}-E_{n+1,n+1})=\td D-n_1+c+x_n\ptl_{x_n}\eqno(2.5)$$
For any $k\in\mbb{Z}$, we denote
$${\msr A}_{\la
k\ra}=\mbox{Span}\:\{x^\al=\prod_{i=1}x_i^{\al_i}\mid\al=(\al_1,...,\al_n)\in\mbb{N}\:^n;\sum_{i=1}^{n_1}\al_i-\sum_{r=n_1+1}^n\al_r=k\}.
\eqno(2.6)$$ Then ${\msr A}=\bigoplus_{k\in\mbb{Z}}{\msr A}_{\la
k\ra}$ and
$${\msr A}_{\la k\ra}=\{f\in{\msr A}\mid
\td D(f)=kf\}.\eqno(2.7)$$ Note that
$${\msr G}_0=\sum_{1\leq i<j\leq
n}(\mbb{F}E_{i,j}+\mbb{F}E_{j,i})+\sum_{r=1}^{n-1}\mbb{F}(E_{r,r}-E_{r+1,r+1})\eqno(2.8)$$
is a Lie subalgebra of $sl(n+1,\mbb{F})$ isomorphic to
$sl(n,\mbb{F})$. The following result was due to Howe \cite{H}.\pse

{\bf Lemma 2.1}. {\it Let $\ell_1,\ell_2\in\mbb{N}$ with $\ell_1>0$,
${\msr A}_{\la -\ell_1\ra}$ is an irreducible highest-weight ${\msr
G}_0$-submodule with highest weight
$\ell_1\lmd_{n_1-1}-(\ell_1+1)\lmd_{n_1}$ and ${\msr A}_{\la
\ell_2\ra}$ is an irreducible highest-weight ${\msr G}_0$-submodule
with highest weight $-(\ell_2+1)\lmd_{n_1}+\ell_2\lmd_{n_1+1}$.}
 \psp

Now we have the first result in this section.\psp

{\bf Theorem 2.2}. {\it The representation $\pi_{c,S}^{\vec 0}$ of
$sl(n+1,\mbb{F})$ is irreducible if any $c\not\in\mbb{Z}$.}

{\it Proof}. Let $k$ be any integer.  For any $0\neq f\in {\msr
A}_{\la k\ra}$, we have
$$0\neq E_{n+1,1}(f)=x_1f\in {\msr A}_{\la k-1\ra}\eqno(2.9)$$
by (2.4), and
$$0\neq E_{n,n+1}(f)=(k+c-n_1)x_nf\in {\msr A}_{\la
k+1\ra}\eqno(2.10)$$ by (2.3). Let ${\msr M}$ be a nonzero
$sl(n+1,\mbb{F})$-submodule of ${\msr A}$. If $k_1,k_2\in\mbb{Z}$
with $k_1\neq k_2$, then the highest weights of ${\msr A}_{\la
k_1\ra}$ and ${\msr A}_{\la k_2\ra}$ are different as ${\msr
G}_0$-modules by Lemma 2.1. So ${\msr A}_{\la k_0\ra}\subset {\msr
M}$ for some $k_0\in\mbb{Z}$. Moreover, (2.9) and (2.10) imply
${\msr A}_{\la k\ra}\subset {\msr M}$ for any $k\in\mbb{Z}$. Hence
${\msr M}={\msr A}$. $\qquad\Box$ \psp

Expressions (2.2)-(2.5) imply the above representation is not of
highest-weight type. Moreover, ${\msr A}$ is a weight
$sl(n+1,\mbb{F})$-module with finite-dimensional weight subspaces.

\psp

{\it Case 2. The representation $\pi_{c,\emptyset}^{\vec a}$ with
$\vec 0\neq \vec a\in\mbb{F}^n$.} \psp

In this case,
$$E_{n+1,i}(fe^{\vec a\cdot\vec
x})=-(\ptl_{x_i}+a_i)(f)e^{\vec a\cdot\vec
x}\qquad\for\;\;i\in\ol{1,n},\;f\in{\msr A}.\eqno(2.11)$$ Thus
$$(E_{n+1,i}+a_i)(fe^{\vec a\cdot\vec
x})=-\ptl_{x_i}(f)e^{\vec a\cdot\vec
x}\qquad\for\;\;i\in\ol{1,n},\;f\in{\msr A}.\eqno(2.12)$$ The second
result in this section.\psp

{\bf Theorem 2.3}. {\it The representation $\pi_{c,\emptyset}^{\vec
a}$ with $\vec 0\neq \vec a\in\mbb{F}^n$ is an irreducible
representation of $sl(n+1,\mbb{F})$ for any $c\in\mbb{F}$.}

{\it Proof}. Let $\msr{A}_k$ be the subspace of homogeneous
polynomials with degree $k$. Set
$$\msr{A}_{\vec a,k}=\msr{A}_ke^{\vec a\cdot\vec
x}\qquad\for\;k\in\mbb{N}.\eqno(2.13)$$ Without loss of generality,
we assume $a_1\neq 0$. Let  ${\msr M}$ be a nonzero
$sl(n+1,\mbb{F})$-submodule of ${\msr A}_{\vec a}$. Take any $0\neq
fe^{\vec a\cdot\vec x}\in \msr{M}$ with $f\in \msr{A}$. By (2.12),
$$\ptl_{x_i}(f)e^{\vec a\cdot\vec
x}\in\msr{M}\qquad\for\;\;i\in\ol{1,n}.\eqno(2.14)$$ By induction,
we have $e^{\vec a\cdot\vec x}\in\msr{M}$; that is, $\msr{A}_{\vec
a,0}\subset\msr{M}$.

Suppose $\msr{A}_{\vec a,\ell}\subset\msr{M}$ for some
$\ell\in\mbb{N}$.  For any $ge^{\vec a\cdot\vec x}\in \msr{A}_{\vec
a,\ell}$,
$$E_{i,1}(ge^{\vec a\cdot\vec x})=x_i(\ptl_{x_1}+a_1)(g)e^{\vec a\cdot\vec x}
=a_1x_ige^{\vec a\cdot\vec x}+x_i\ptl_{x_1}(g)e^{\vec a\cdot\vec
x}\in\msr{M}\qquad\for\;\;i\in\ol{2,n}\eqno(2.15)$$ by (1.4). Since
$x_i\ptl_{x_1}(g)e^{\vec a\cdot\vec x}\in \msr{A}_{\vec
a,\ell}\subset\msr{M}$, we have
$$x_ige^{\vec a\cdot\vec x}\in\msr{M}\qquad\for\;\;i\in\ol{2,n}.\eqno(2.16)$$

On the other hand,
$$(E_{1,1}-E_{2,2})(ge^{\vec a\cdot\vec x})=a_1x_1ge^{\vec a\cdot\vec
x}+(x_1\ptl_{x_1}-x_2\ptl_{x_2}-a_2x_2)(g)e^{\vec a\cdot\vec
x}\in\msr{M}\eqno(2.17)$$ by (1.5).  Our assumption says that
$(x_1\ptl_{x_1}-x_2\ptl_{x_2})(g)e^{\vec a\cdot\vec x}\in
\msr{A}_{\vec a,\ell}\subset\msr{M}$. According to (2.16),
$-a_2x_2(g)e^{\vec a\cdot\vec x}\in\msr{M}$. Therefore,
$$x_1ge^{\vec a\cdot\vec
x}\in\msr{M}.\eqno(2.18)$$ Expressions (2.17) and (2.18) imply
$\msr{A}_{\vec a,\ell+1}\subset\msr{M}$. By induction,
$\msr{A}_{\vec a,\ell}\subset\msr{M}$ for any $\ell\in\mbb{N}$. So
$\msr{A}_{\vec a}=\msr{M}$. Hence $\msr{A}_{\vec a}$ is an
irreducible $sl(n+1,\mbb{F})$-module. $\qquad\Box$\psp

{\it Case 3. The representation $\pi_{c,S}^{\vec a}$ with  $a_i\neq
0$ for some $i\in \ol{1,n}\setminus S$ or $|S|>1$.}\psp

The following is the third result in this section.\psp

{\bf Theorem 2.4}. {\it Under the above assumption, the
representation $\pi_{c,S}^{\vec a}$ with $\vec 0\neq \vec
a\in\mbb{F}^n$ is an irreducible representation of
$sl(n+1,\mbb{F})$.}

{\it Proof}. Without loss of generality, we assume $S=\ol{1,n_1}$
for some $n_1\in\ol{1,n_1}$ and $n_1<n$. Let ${\msr M}$ be a nonzero
$sl(n+1,\mbb{F})$-submodule of ${\msr A}_{\vec a}$. By (2.4) and
(2.11)-(2.14), there exists $0\neq fe^{\vec a\cdot\vec x}\in\msr{M}$
with $f\in\mbb{F}[x_1,...,x_{n_1}]$. \pse

{\it Subcase (1). $a_i\neq 0$ for some $i\in\ol{n_1+1,n}$.}\pse

By symmetry, we can assume $a_n\neq 0$. According to (2.2),
$$E_{i,n}(fe^{\vec a\cdot\vec
x})=(\ptl_{x_i}+a_i)(\ptl_{x_n}+a_n)(f)e^{\vec a\cdot\vec
x}=a_ia_nfe^{\vec a\cdot\vec x}+a_n\ptl_{x_i}(f)e^{\vec a\cdot\vec
x}\;\;\for\;i\in\ol{1,n_1}.\eqno(2.19)$$ Thus
$$(a_n^{-1}E_{i,n}-a_i)(fe^{\vec a\cdot\vec
x})=\ptl_{x_i}(f)e^{\vec a\cdot\vec
x}\in\msr{M}\;\;\for\;i\in\ol{1,n_1}.\eqno(2.20)$$ By induction on
the degree of $f$, we get $e^{\vec a\cdot\vec x}\in\msr{M}$; that
is, $\msr{A}_{\vec a,0}\subset\msr{M}$.

 The arguments in (2.15)-(2.18) yield
$$\mbb{F}[x_{n_1+1},...,x_n]e^{\vec a\cdot\vec
x}\subset\msr{M}.\eqno(2.21)$$ According to (2.4),
$$E_{n+1,1}^{\ell_1}\cdots E_{n+1,n_1}^{\ell_{n_1}}(\mbb{F}[x_{n_1+1},...,x_n]e^{\vec a\cdot\vec
x})=x_1^{\ell_1}\cdots
x_{n_1}^{\ell_{n_1}}(\mbb{F}[x_{n_1+1},...,x_n]e^{\vec a\cdot\vec
x})\subset\msr{M}\eqno(2.22)$$ for $\ell_i\in\mbb{N}$ with
$i\in\ol{1,n_1}.$ Thus ${\msr A}_{\vec a}=\msr{M}$. So ${\msr
A}_{\vec a}$ is $\msr{A}_{\vec a}$ is an irreducible
$sl(n+1,\mbb{F})$-module. \pse

{\it Subcase (2). $a_i= 0$ for any $i\in\ol{n_1+1,n}$ and
$n_1>1$.}\pse

By the transformation
$$\vec x\mapsto T\vec x ,\;\;A\mapsto TAT^{-1}\eqno(2.23)$$
with $A\in sl(n+1,\mbb{F})$ for some $n\times n$ orthogonal matrix
$T$, we can assume $a_1\neq 0$ and $a_i=0$ for $i\in\ol{2,n}$. Note
that
$$E_{1,2}(fe^{\vec a\cdot\vec x})=-x_2(\ptl_{x_1}+a_1)(f)e^{\vec a\cdot\vec
x}\in\msr{M}\eqno(2.24)$$ by (2.2) and
$$E_{n+1,2}(fe^{\vec a\cdot\vec x})=x_2fe^{\vec a\cdot\vec x}\in\msr{M}\eqno(2.25)$$ by
(2.4). Thus
$$x_2\ptl_{x_1}(f)e^{\vec a\cdot\vec
x}\in\msr{M}.\eqno(2.26)$$ Repeatedly applying (2.26) if necessary,
we can assume $f\in\mbb{F}[x_2,....,x_n]$. We apply the arguments in
the proof of Theorem 2.2 to the Lie subalgebra
$$\msr{L}=\sum_{2\leq i<j\leq
n+1}(\mbb{F}E_{i,j}+\mbb{F}E_{j,i})+\sum_{r=2}^n\mbb{F}(E_{r,r}-E_{r+1,r+1})\eqno(2.27)$$
and obtain
$$\mbb{F}[x_2,....,x_n]e^{\vec a\cdot\vec x}\subset \msr{M}.\eqno(2.28)$$
According to (2.4),
$$E_{n+1,1}^\ell(\mbb{F}[x_2,....,x_n]e^{\vec a\cdot\vec x})=x_1^\ell(\mbb{F}[x_2,....,x_n]e^{\vec a\cdot\vec x})
\subset \msr{M}\eqno(2.29)$$ for any $\ell\in\mbb{N}$. Therefore
$\msr{M}=\msr{A}_{\vec a}$. So $\msr{A}_{\vec a}$ is an irreducible
$sl(n+1,\mbb{F})$-module. $\qquad\Box$\psp

With the representation $\pi_{c,S}^{\vec 0}$, ${\msr A}$ is an
infinite-dimensional  weight $sl(n+1,\mbb{F})$-module with
finite-dimensional weight subspaces by (2.2) and (2.5). Now Theorem
1 follows from Theorems 2.2, 2.3 and 2.4.

\section{Proof of Theorem 2}

$\quad\;$ Assume that $n=2m+1>1$ is an odd integer. In this section,
we will give the proof of Theorem 2.

 Recall that the symplectic
Lie algebras \begin{eqnarray*} \qquad sp(2m+2,\mbb{F})&=&\sum_{1\leq
r\leq s\leq
m+1}[\mbb{F}(E_{r,m+1+s}+E_{s,m+1+r})+\mbb F(E_{m+1+r,s}+E_{m+1+s,r})]\\
&&+
\sum_{i,j=1}^{m+1}\mbb{F}(E_{i,j}-E_{m+1+j,m+1+i}).\hspace{5.3cm}(3.1)\end{eqnarray*}
For convenience, we rednote
$$x_0=x_{m+1},\;\;y_i=x_{m+1+i}\qquad\for\;i\in\ol{1,m}.\eqno(3.2)$$
In particular,
$$D=\sum_{i=0}^mx_i\ptl_{x_i}+\sum_{r=1}^my_r\ptl_{y_r}.\eqno(3.3)$$
According to (1.4) and (1.5), we have the representation $\pi_c$ of
$sp(2m+2,\mbb{F})$:
$$\pi_c(E_{i,j}-E_{m+1+j,m+1+i})=x_i\ptl_{x_j}-y_j\ptl_{y_i},\;\;\pi_c(E_{i,m+1+j}+E_{j,m+1+i})
=x_i\ptl_{y_j}+x_j\ptl_{y_i},\eqno(3.4)$$
$$\pi_c(E_{2m+2,m+1})=-\ptl_{x_0},\;\;\pi_c(E_{m+1,i}-E_{m+1+i,2m+2})=x_0\ptl_{x_i}-y_i(D+c),
\eqno(3.5)$$
$$\pi_c(E_{i,m+1}-E_{2m+2,m+1+i})=x_i\ptl_{x_0}+\ptl_{y_i},\;\;\pi_c(E_{2m+2,i}+E_{m+1+i,m+1})=y_i\ptl_{x_0}-\ptl_{x_i},\eqno(3.6)$$
$$\pi_c(E_{m+1+i,j}+E_{m+1+j,i})=y_i\ptl_{x_j}+y_j\ptl_{x_i},\;\;\pi_c(E_{m+1,m+1}-E_{2m+2,2m+2})=D+x_0\ptl_{x_0}+c,
\eqno(3.7)$$
$$\pi_c(E_{m+1,m+1+i}+E_{i,2m+2})=x_0\ptl_{y_i}+x_i(D+c),\;\;\pi_c(E_{i,i}-E_{m+i,m+i})
=x_i\ptl_{x_i}-y_i\ptl_{y_i},\eqno(3.8)$$
$$\pi_c(E_{m+1,2m+2})=x_0(D+c)\eqno(3.9)$$
 for
$i,j\in\ol{1,m}$.

Denote \begin{eqnarray*} \qquad\qquad {\msr K}&=&\sum_{1\leq r\leq
s\leq
m}[\mbb{F}(E_{r,m+1+s}+E_{s,m+1+r})+\mbb{F}(E_{m+1+r,s}+E_{m+1+s,r})]\\
&&+
\sum_{i,j=1}^m\mbb{F}(E_{i,j}-E_{m+1+j,m+1+i}),\hspace{6.3cm}(3.10)\end{eqnarray*}
which is a Lie subalgebra of $sp(2m+2,\mbb{F})$ isomorphic to
$sp(2m,\mbb{F})$. We will prove Theorem 2 case by case.\psp

{\it Case 1. $\vec a=\vec 0$ and $S=\emptyset$.}\psp

 Let $\msr{B}=\mbb{F}[x_1,...,x_n,y_1,...,y_n]$. Denote by
 $\msr{B}_k$ the subspace of homogeneous polynomials with degree
 $k$. First we have the following  well-known result
(e.g., cf. \cite{FH}).\pse

{\bf Lemma 3.1}. {\it For any $k\in\mbb{N}$, ${\msr B}_k$ forms a
finite-dimensional irreducible $\msr{K}$-module with highest weight
$k\lmd_1$.}\pse

Set
$${\msr A}_{(\ell)}=\sum_{i=0}^\ell {\msr
A}_i\qquad\for\;\;\ell\in\mbb{N}.\eqno(3.11)$$ Take the Cartan
subalgebra
$$H=\sum_{i=1}^{m+1}\mbb{F}(E_{i,i}-E_{m+1+i,m+1+i})\eqno(3.12)$$ of
$sp(2m+2,\mbb{F})$. Define $\{\ves_1,...,\ves_{n+1}\}\subset H^\ast$
by:
$$\ves_j(E_{i,i}-E_{n+1+i,n+1+i})=\dlt_{i,j}.\eqno(3.13)$$
 Recall that
the representation $\pi_{c,\emptyset}^{\vec 0}$ of
$sp(2m+2,\mbb{F})$ is the representation $\pi_c$ (cf. (3.4)-(3.9))
on the space $\msr{A}=\mbb{F}[x_0,...,x_m,y_1,...,y_m]$. Then we
have: \psp

{\bf Theorem 3.2}. {\it If $c\not\in -\mbb{N}$, the representation
$\pi_{c,\emptyset}^{\vec 0}$ of $sp(2m+2,\mbb{F})$ given in
(3.4)-(3.9) is a highest-weight irreducible representation with
highest-weight
 $-c\lmd_1$. When $-c=\ell\in\mbb{N}$, ${\msr A}_{(\ell)}$ is
 a finite-dimensional irreducible
$sp(2m+2,\mbb{F})$-module with highest weight $\ell\lmd_n$ and
${\msr A}/{\msr A}_{(\ell)}$ is an irreducible highest weight
$sp(2m+2,\mbb F)$-module with highest weight
$-(\ell+2)\lmd_1+(\ell+1)\lmd_2$, where $\lmd_i$ is the $i$th
fundamental weight of $sp(2m+2,\mbb{F})$.}

{\it Proof}.  Observe that
$${\msr A}_k=\sum_{s=0}^kx_0^s{\cal
B}_{k-s}\qquad\for\;\;k\in\mbb{N}.\eqno(3.14)$$ Let $\msr{M}$ be a
nonzero $sp(2m+2,\mbb{F})$-submodule of ${\msr A}$. Take any $0\neq
f\in\msr{M}$. Repeatedly applying the first equation in (3.5) and
(3.6) to $f$, we obtain $1\in \msr{M}$. Note
$$(E_{m+1,2m+2})^k(1)=[\prod_{r=0}^{k-1}(r+c)]x_0^k\in
\msr{M}\eqno(3.15)$$ by (3.9) and
$$(E_{1,m+1}-E_{2m+2,m+1+1})^s(x_0^k)=[\prod_{i=0}^{s-1}(k-i)x_0^{k-s}x_1^s\qquad\for\;\;s\in\ol{1,k}\eqno(3.16)$$
by the first equation in (3.6). Suppose $c\not\in-\mbb{N}$. Then
(3.17) yields
$$x_0^k\in \msr{M}\qquad\for\;\;k\in\mbb{N}.\eqno(3.17)$$
Moreover, (3.16) with $k=r+s$ gives
$$x_0^rx_1^s\in V\qquad\for\;\;r,s\in\mbb{N}.\eqno(3.18)$$
Furthermore,
$$U({\cal K})(x_0^rx_1^s)=x_0^r{\msr B}_s\subset \msr{M}\eqno(3.19)$$
by Lemma 3.1. Thus
$${\msr A}=\sum_{r,s=0}^\infty x_0^r{\msr B}_s\subset
\msr{M};\eqno(3.20)$$ that is, $\msr{M}={\msr A}$. So ${\msr A}$ is
an irreducible $sp(2m+2,\mbb F)$-module and $1$ is its
highest-weight vector with weight $-c\lmd_1$ with respect to the
following simple positive roots
$$\{\ves_{n+1}-\ves_n,\ves_n-\ves_{n-1},...,\ves_2-\ves_1,2\ves_1
\}.\eqno(3.21)$$

Next we assume $c=-\ell$ with $\ell\in\mbb{N}$. Since
$$E_{m+1,2m+2}|_{{\msr A}_\ell}=0,\;\;
(E_{m+1,i}-E_{m+1+i,2m+2})|_{{\msr
A}_\ell}=x_0\ptl_{x_i},\eqno(3.22)$$
$$(E_{m+1,m+1+i}+E_{i,2m+2})|_{{\msr A}_\ell}=x_0\ptl_{y_i}\eqno(3.23)$$ by (3.5), (3.8) and (3.9),
${\msr A}_{(\ell)}$ is
 a finite-dimensional $sp(2m+2,\mbb{F})$-module.
Let $\msr{M}$ be a nonzero $sp(2m+2,\mbb{F})$-submodule of ${\msr
A}_{(\ell)}$. By (3.15),
$$x_0^k\in \msr{M}\qquad\for\;\;k\in\ol{0,\ell}.\eqno(3.24)$$
Moreover, (3.16) with $k=r+s$ gives
$$x_0^rx_1^s\in \msr{M} \qquad\for\;\;r,s\in\ol{0,\ell}\;\;\mbox{such that}\;\;r+s\leq \ell.\eqno(3.25)$$
 Thus
$$\msr{A}_{(\ell)}=\sum_{r=0}^\ell\sum_{s=0}^{\ell-r} x_0^r{\cal B}_s\subset
\msr{M}\eqno(3.26)$$ by Lemma 3.1; that is,
$\msr{M}=\msr{A}_{(\ell)}$. So $\msr{A}_{(\ell)}$ is an irreducible
$sp(2n+2,\mbb{F})$-module and $1$ is again its highest-weight
vector.

Consider the quotient $sp(2n+2,\mbb{F})$-module
$\msr{A}/\msr{A}_{(\ell)}$. Let $W\supset \msr{A}_{(\ell)}$ be an
$sp(2n+2,\mbb{F})$-submodule of $\msr{A}$ such that $W\neq
\msr{A}_{(\ell)}$. Take any $f\in W\setminus \msr{A}_{(\ell)}$.
Repeatedly applying (3.6) and the first equation in (3.5) to $f$ if
necessary, we can assume $f\in {\cal B}_{\ell+1}$. Since ${\cal
B}_{\ell+1}$ is an irreducible $\msr{K}$-module, we have
$${\cal B}_{\ell+1}\subset W.\eqno(3.27)$$ In particular,
$x_1^{\ell+1}\in W$. According to (3.8),
$$(E_{m+1,m+2}+E_{1,2m+2})^r(x_1^{\ell+1})=r!x_1^{\ell+1+r}\in
W\qquad \for\;\;0<r\in\mbb{Z}.\eqno(3.28)$$ Since ${\cal
B}_{\ell+1+r}\ni x_1^{\ell+1+r}$ is an irreducible $\msr{K}$-module,
we have
$${\cal B}_{\ell+1+r}\subset W.\eqno(3.29)$$

Suppose that
$$x_0^r{\cal B}_s \subset
W\qquad\for\;\;r\in\ol{0,k}\;\mbox{and}\;s\in\mbb{N}\;\mbox{such
that}\;r+s\geq\ell+1.\eqno(3.30)$$ Fix such $r$ and $s$. Observe
$x_0^rx_1^{s-1}y_1\in x_0^r{\cal B}_s\subset W$. Using the first
equation in (3.8), we get
$$(E_{m+1,m+2}+E_{1,2m+2})(x_0^rx_1^{s-1}y_1)=(r+s-\ell)x_0^rx_1^sy_1
+x_0^{r+1}x_1^{s-1}\in W.\eqno(3.31)$$ By the assumption (3.30),
$(r+s-\ell)x_0^rx_1^sy_1\in x_0^r{\cal B}_{s+1}\subset W$. So
$$x_0^{r+1}x_1^{s-1}\in W\bigcap x_0^{r+1}{\cal
B}_{s-1}.\eqno(3.32)$$ Since $ x_0^{r+1}{\cal B}_{s-1}$ is an
irreducible $\msr{K}$-module, we get
$$ x_0^{r+1}{\cal B}_{s-1}\subset W.\eqno(3.33)$$
By induction on $r$, we prove
$$x_0^r{\cal B}_s \subset
W\qquad\for\;\;r,s\in\mbb{N}\;\mbox{such that}\;r+s\geq
\ell+1.\eqno(3.34)$$ According to (3.14),
$$\sum_{k=\ell+1}^\infty\msr{A}_k\subset W.\eqno(3.35)$$
Since $W\supset \msr{A}_{(\ell)}$, we have $W=\msr{A}$. So
$\msr{A}/\msr{A}_{(\ell)}$ is an irreducible
$sp(2m+2,\mbb{F})$-module. Moreover, $x_n^{\ell+1}$ is a highest
weight vector of weight $-(\ell+2)\lmd_1+(\ell+1)\lmd_2$ with
respect to (3.21). $\qquad\Box$\psp

{\it Case 2. $\vec a\neq\vec0,\;a_{m+1}=0$ and $S=\emptyset$.}\psp

For simplicity, we redenote
$$
b_i=a_{m+1+i}\qquad\for\;\;i\in\ol{1,m}.\eqno(3.36)$$ Recall that
the representation $\pi_{c,\emptyset}^{\vec a}$ of
$sp(2m+2,\mbb{F})$ is the representation $\pi_c$ (cf. (3.4)-(3.9))
on the space $\msr{A}_{\vec a}$ (cf. (1.7)). Our second result in
this section is:\psp

{\bf Theorem 3.3}. {\it The representation $\pi_{c,\emptyset}^{\vec
a}$ with $\vec 0\neq \vec a\in\mbb{F}^n$ and $a_{m+1}=0$ is an
irreducible representation of $sp(2m+2,\mbb{F})$ for any
$c\in\mbb{F}$.}

{\it Proof}. By symmetry, we may assume $a_1\neq 0$. Let ${\msr M}$
be a nonzero $sp(2m+2,\mbb{F})$-submodule of ${\msr A}_{\vec a}$.
Take any $0\neq fe^{\vec a\cdot\vec x}\in \msr{M}$ with $f\in
\msr{A}$. By the assumption $a_0=0$, (3.5) and (3.6),
$$E_{2m+2,m+1}(fe^{\vec a\cdot\vec x})=-\ptl_{x_0}(f)e^{\vec a\cdot\vec
x} \in \msr{M},\eqno(3.37)$$
$$(E_{i,m+1}-E_{2m+2,m+1+i}-b_i)(fe^{\vec a\cdot\vec
x})=[\ptl_{y_i}(f)+x_i\ptl_{x_0}(f)]e^{\vec a\cdot\vec x} \in
\msr{M},\eqno(3.38)$$
$$(E_{2m+2,i}+E_{m+1+i,m+1}+a_i)(fe^{\vec a\cdot\vec
x})=[-\ptl_{x_i}(f)+y_i\ptl_{x_0}(f)]e^{\vec a\cdot\vec x} \in
\msr{M}\eqno(3.39)$$ for $i\in\ol{1,m}$. Repeatedly applying
(3.37)-(3.39), we obtain $e^{\vec a\cdot\vec x}\in \msr{M}$.
Equivalently, $\msr{A}_{\vec a,0}\subset\msr{M}$ (cf. (2.13)).

Suppose $\msr{A}_{\vec a,\ell}\subset\msr{M}$ for some
$\ell\in\mbb{N}$.  For any $ge^{\vec a\cdot\vec x}\in \msr{A}_{\vec
a,\ell}$,
$$(E_{i,1}-E_{m+2,m+1+i})(ge^{\vec a\cdot\vec x})=[a_1x_i-b_iy_1+x_i\ptl_{x_1}-y_1\ptl_{x_i}](g)e^{\vec a\cdot\vec x}
\in\msr{M}\eqno(3.40)$$ by (3.4) and
$$(E_{m+1+i,1}+E_{m+2,i})(ge^{\vec a\cdot\vec x})=[a_1y_i+a_iy_1+y_i\ptl_{x_1}+y_1\ptl_{x_i}](g)e^{\vec a\cdot\vec x}
\in\msr{M}\eqno(3.41)$$ by the first equation in (3.7), where
$i\in\ol{1,m}$. Since
$$(x_i\ptl_{x_1}-y_1\ptl_{x_i})(g)e^{\vec a\cdot\vec x},\;(y_i\ptl_{x_1}+y_1\ptl_{x_i})(g)e^{\vec a\cdot\vec
x}\in \msr{A}_{\vec a,\ell}\subset\msr{M},\eqno(3.42)$$ we have
$$(a_1x_i-b_iy_1)ge^{\vec a\cdot\vec x},\;(a_1y_i+a_iy_1)ge^{\vec a\cdot\vec x}
\in\msr{M}\eqno(3.43)$$ for $i\in\ol{1,m}$. The above second
equation with $i=1$ gives
$$2a_1y_1ge^{\vec a\cdot\vec x}
\in\msr{M}\Rightarrow y_1ge^{\vec a\cdot\vec x}
\in\msr{M}.\eqno(3.44)$$ Thus (3.43) yields
$$x_ige^{\vec a\cdot\vec x},y_ige^{\vec a\cdot\vec x}
\in\msr{M}\qquad\for\;i\in\ol{1,m}.\eqno(3.45)$$

According to the second equation in (3.5),
\begin{eqnarray*}& &(E_{m+1,1}-E_{m+2,2m+2})(ge^{\vec a\cdot\vec
x})\\
&=&[a_1x_0-\sum_{i=1}^m(a_ix_i+b_iy_i)y_1+x_0\ptl_{x_1}-y_1(D+c)](g)e^{\vec
a\cdot\vec x}\in\msr{M}.\hspace{3.2cm}(3.46)\end{eqnarray*}
Replacing $ge^{\vec a\cdot\vec x}\in \msr{A}_{\vec a,\ell}$ by
$ge^{\vec a\cdot\vec x}\in \sum_{i=1}^m(x_i\msr{A}_{\vec
a,\ell}+y_i\msr{A}_{\vec a,\ell})$ in (3.40)-(3.45), we obtain
$$x_iy_1ge^{\vec a\cdot\vec x},y_iy_1ge^{\vec a\cdot\vec x}
\in\msr{M}\qquad\for\;i\in\ol{1,m}.\eqno(3.47)$$ Since $D(g)=\ell g$
and $x_0\ptl_{x_1}(g)e^{\vec a\cdot\vec x}\in \msr{A}_{\vec
a,\ell}\subset\msr{M}$, we have
$$[-\sum_{i=1}^m(a_ix_i+b_iy_i)y_1+x_0\ptl_{x_1}-y_1(D+c)](g)e^{\vec
a\cdot\vec x}\in\msr{M}.\eqno(3.48)$$ Hence (3.46) yields
$x_0ge^{\vec a\cdot\vec x}\in\msr{M}$. Therefore, $\msr{A}_{\vec
a,\ell+1}\subset\msr{M}$. By induction, $\msr{A}_{\vec
a,\ell}\subset\msr{M}$ for any $\ell\in\mbb{N}$. So $\msr{A}_{\vec
a}=\msr{M}$. Hence $\msr{A}_{\vec a}$ is an irreducible
$sp(2m+2,\mbb{F})$-module. $\qquad\Box$\psp

{\it Case 3. $\vec a=\vec 0$ and $S\neq\emptyset$.}\psp

By symmetry and the assumption (1.8), we can assume
$$S=\ol{1,m_1}\bigcup\ol{m_2+1,m}, \qquad
m_1,m_2\in\ol{1,m}\;\mbox{and}\;m_1\leq m_2,\eqno(3.49)$$ where we
treat $\ol{m+1,m}=\emptyset$ when $m_2=m$. Set
 $$\td D=x_0\ptl_{x_0}+\sum_{r=m_1+1}^mx_r\ptl_{x_r}-\sum_{i=1}^{m_1}x_i\ptl_{x_i}+
\sum_{i=1}^{m_2}y_i\ptl_{y_i}-\sum_{r=m_2+1}^my_r\ptl_{y_r}\eqno(3.50)$$
and
$$\td c=c+m_2-m_1-m.\eqno(3.51)$$
Then we have the following representation $\pi_{c,S}$ of the Lie
algebra $sp(2m+2,\mbb{F})$  determined by
$$\pi_{c,S}(E_{i,j}-E_{m+1+j,m+1+i})=E_{i,j}^x-E_{j,i}^y\eqno(3.52)$$ with
$$E_{i,j}^x=\left\{\begin{array}{ll}-x_j\ptl_{x_i}-\delta_{i,j}&\mbox{if}\;
i,j\in\ol{1,m_1};\\ \ptl_{x_i}\ptl_{x_j}&\mbox{if}\;i\in\ol{1,m_1},\;j\in\ol{m_1+1,m};\\
-x_ix_j &\mbox{if}\;i\in\ol{m_1+1,m},\;j\in\ol{1,m_1};\\
x_i\partial_{x_j}&\mbox{if}\;i,j\in\ol{m_1+1,m}
\end{array}\right.\eqno(3.53)$$
and
$$E_{i,j}^y=\left\{\begin{array}{ll}y_i\ptl_{y_j}&\mbox{if}\;
i,j\in\ol{1,m_2};\\ -y_iy_j&\mbox{if}\;i\in\ol{1,m_2},\;j\in\ol{m_2+1,m};\\
\ptl_{y_i}\ptl_{y_j} &\mbox{if}\;i\in\ol{m_2+1,m},\;j\in\ol{1,m_2};\\
-y_j\partial_{y_i}-\delta_{i,j}&\mbox{if}\;i,j\in\ol{m_2+1,m},
\end{array}\right.\eqno(3.54)$$
and
$$\pi_{c,S}(E_{i,m+1+j})=\left\{\begin{array}{ll}
\ptl_{x_i}\ptl_{y_j}&\mbox{if}\;i\in\ol{1,m_1},\;j\in\ol{1,m_2},\\
-y_j\ptl_{x_i}&\mbox{if}\;i\in\ol{1,m_1},\;j\in\ol{m_2+1,m},\\
x_i\ptl_{y_j}&\mbox{if}\;i\in\ol{m_1+1,m},\;j\in\ol{1,m_2},\\
-x_iy_j&\mbox{if}\;i\in\ol{m_1+1,m},\;j\in\ol{m_2+1,m},\end{array}\right.\eqno(3.55)$$
$$\pi_{c,S}(E_{m+1+i,j})=\left\{\begin{array}{ll}
-x_jy_i&\mbox{if}\;j\in\ol{1,m_1},\;i\in\ol{1,m_2},\\
-x_j\ptl_{y_i}&\mbox{if}\;j\in\ol{1,m_1},\;i\in\ol{m_2+1,m},\\
y_i\ptl_{x_j}&\mbox{if}\;j\in\ol{m_1+1,m},\;i\in\ol{1,m_2},\\
\ptl_{x_j}\ptl_{y_i}&\mbox{if}\;j\in\ol{m_1+1,m},\;i\in\ol{m_2+1,m},\end{array}\right.\eqno(3.56)$$
 $$\pi_{c,S}(E_{2m+2,m+1})=-\ptl_{x_0},\;\;\pi_{c,S}(E_{2m+2,m+1})=x_0(\td
D+\td c),\eqno(3.57)$$
$$\pi_{c,S}(E_{i,m+1}-E_{2m+2,m+1+i})=\left\{\begin{array}{ll}\ptl_{x_0}\ptl_{x_i}+\ptl_{y_i}&\mbox{if}\;\;i\in\ol{1,m_1},
\\ x_i\ptl_{x_0}+\ptl_{y_i}&\mbox{if}\;\;i\in\ol{m_1+1,m_2},\\ x_i\ptl_{x_0}-y_i&
\mbox{if}\;\;i\in\ol{m_2+1,m},\end{array}\right.\eqno(3.58)$$
$$\pi_{c,S}(E_{2m+2,i}+E_{m+1+i,m+1})=\left\{\begin{array}{ll}y_i\ptl_{x_0}+x_i&\mbox{if}\;\;i\in\ol{1,m_1},
\\ y_i\ptl_{x_0}-\ptl_{x_i}&\mbox{if}\;\;i\in\ol{m_1+1,m_2},\\ \ptl_{x_0}\ptl_{y_i}-\ptl_{x_i}&
\mbox{if}\;\;i\in\ol{m_2+1,m},\end{array}\right.\eqno(3.59)$$
$$\pi_{c,S}(E_{m+1,i}-E_{m+1+i,2m+2})=\left\{\begin{array}{ll}-x_0x_i-y_i(\td D+\td c)
&\mbox{if}\;\;i\in\ol{1,m_1},
\\ x_0\ptl_{x_i}-y_i(\td D+\td c)&\mbox{if}\;\;i\in\ol{m_1+1,m_2},\\ x_0\ptl_{x_i}-(\td D+\td
c-1)\ptl_{y_i} &
\mbox{if}\;\;i\in\ol{m_2+1,m},\end{array}\right.\eqno(3.60)$$
$$\pi_{c,S}(E_{m+1,m+1+i}+E_{i,2m+2})=\left\{\begin{array}{ll}x_0\ptl_{y_i}+(\td D+\td
c-1)\ptl_{x_i} &\mbox{if}\;\;i\in\ol{1,m_1},
\\ x_0\ptl_{y_i}+x_i(\td D+\td c)&\mbox{if}\;\;i\in\ol{m_1+1,m_2},\\ -x_0y_i+x_i(\td D+\td c)
& \mbox{if}\;\;i\in\ol{m_2+1,m},\end{array}\right.\eqno(3.61)$$
$$\pi_{c,S}(E_{m+1,m+1}-E_{2m+2,2m+2})=\td D+x_0\ptl_{x_0}+\td c,
\eqno(3.62)$$ for $i,j\in\ol{1,m}$.

Recall $\msr{B}=\mbb{F}[x_1,...,x_n,y_1,...,y_n]$. Set
$${\msr B}_{\la k\ra}={\msr A}_{\la
k\ra}\bigcap \msr{B}\qquad\for\;\;k\in\mbb{Z}\eqno(3.63)$$ (cf.
(2.6)). Then $\msr{B}=\bigoplus_{k\in\mbb{Z}}\msr B_{\la k\ra}$ is a
$\mbb{Z}$-graded space. The following result is due to
\cite{LX2}:\psp

{\bf Lemma 3.4}. {\it Assume $m\geq 2$. Let $k\in\mbb{Z}$. If
$m_1<m_2$ or $k\neq 0$, the subspace ${\msr B}_{\la k\ra}$  is an
irreducible ${\msr K}$-submodule (cf. (3.10)). When $m_1=m_2$, the
subspace ${\msr B}_{\la 0\ra}$ is a direct sum of two irreducible
$\msr{K}$-submodules.} \psp

In fact, any pair of  the irreducible submodules in the above are
not isomorphic $\msr{K}$-modules because they have distinct weight
sets of singular vectors with respect to the Lie subalgebra
$\sum_{i,j=1}^m\mbb{F}(E_{i,j}-E_{m+1+j,m+1+i})\cong gl(m,\mbb{F})$
(cf. \cite{LX1}). When $m=m_1=m_2=1$, $\msr{B}=\mbb{F}[x_1,x_2]$ and
$$\pi_{c,S}(\msr{K})=\mbb{F}(x_1\ptl_{x_1}+y_1\ptl_{y_1}+1)+\mbb{F}x_1y_1+\mbb{F}\ptl_{x_1}\ptl_{y_1}.
\eqno(3.64)$$ So all ${\msr B}_{\la k\ra}$ with $k\in\mbb{Z}$ are
irreducible ${\msr K}$-submodules. Recall the representation
$\pi_{c,S}^{\vec 0}$ of $sp(2m+2,\mbb{F})$ is the representation
$\pi_{c,S}$ (cf. (3.52)-(3.62)) on $\msr A$. The following is the
third result in this section.\psp

{\bf Theorem 3.5}. {\it The representation $\pi_{c,S}^{\vec 0}$ of
$sp(2m+2,\mbb{F})$ is irreducible if $c\not\in\mbb{Z}$.}

{\it Proof}. Let $\msr{M}$ be any nonzero
$sp(2m+2,\mbb{F})$-submodule of $\msr{A}$. Repeatedly applying
$E_{2m+2,m+1}$ to $\msr{M}$ by the first equation in (3.57), we get
$$\msr{A}\bigcap {\msr B}\neq\{0\}.\eqno(3.65)$$
According to (3.62),
$${\msr B}_{\la k\ra}=\{f\in\msr B\mid
(E_{m+1,m+1}-E_{2m+2,2m+2})(f)=(k+\td c)f\}.\eqno(3.66)$$ Thus
$$\msr M=\bigoplus_{k\in\mbb{Z}}\msr M\bigcap {\msr B}_{\la
k\ra}.\eqno(3.67)$$ If $\msr M\bigcap\msr B_{\la 0\ra}\neq \{0\}$,
then (3.59) gives
$$(E_{2m+2,1}+E_{m+2,m+1})(\msr M\bigcap\msr B_{\la 0\ra})=
x_1(\msr M\bigcap\msr B_{\la 0\ra})\subset \msr M\bigcap\msr B_{\la
-1\ra}.\eqno(3.68)$$ Thus we always have $\msr M\bigcap\msr B_{\la
k\ra}\neq\{0\}$ for some $0\neq k\in\mbb{Z}$. According to Lemma 3.4
and (3.64), $\msr B_{\la k\ra}$ is an irreducible $\msr K$-module.
So
$$\msr B_{\la k\ra}\subset\msr M.\eqno(3.69)$$

Next (3.58) yields
$$\msr B_{\la k-r\ra}=(\ptl_{y_1})^r(\msr  B_{\la
k\ra})=(E_{1,m+1}-E_{2m+2,m+2})^r(\msr B_{\la k\ra})\subset \msr
M\qquad\for\;\;r\in\mbb{N}.\eqno(3.70)$$ On the other hand,  if
 $\msr B_{\la \ell\ra}\subset \msr
M$, then the assumption $c\not\in\mbb Z$ and the second equation in
(3.57) give
$$x_0^r\msr B_{\la \ell\ra}=(E_{2m+2,m+1})^r(\msr B_{\la
\ell\ra})\subset \msr M \qquad\for\;\;r\in\mbb{N}.\eqno(3.71)$$
Suppose that for some $s\in\mbb{Z}$,
$$x_0^r\msr B_{\la s\ra},x_0^r{\msr B}_{\la s-1\ra}
\subset\msr M \qquad\for\;\;r\in\mbb{N}.\eqno(3.72)$$ For any
$\ell\in\mbb{N}$,
$$x_0^\ell{\msr B}_{\la s+1\ra}=(\td D+\td
c-1)\ptl_{x_1}(x_0^\ell{\msr B}_{\la
s\ra})=[E_{m+1,m+2}+E_{1,2m+2}-x_0\ptl_{y_1}](x_0^\ell{\msr B}_{\la
s\ra})\eqno(3.73)$$ by (3.61). Note
$$x_0\ptl_{y_1}(x_0^\ell{\msr
B}_{\la s\ra})=x_0^{\ell+1}{\msr B}_{\la s-1\ra}\subset \msr
M.\eqno(3.74)$$ Thus (3.73) leads to
$$x_0^\ell{\msr B}_{\la s+1\ra}\subset\msr M.\eqno(3.75)$$

By (3.70)-(3.75) and induction on $s$, we prove
$$x_0^r{\msr B}_{\la k\ra}\subset
W\qquad\for\;\;x_0\in\mbb{N},\;k\in\mbb{Z}.\eqno(3.76)$$ So $\msr
M=\msr A$. Therefore, $\msr A$ is an irreducible
$sp(2m+2,\mbb{F})$-module. $\qquad\Box$\psp

{\bf Remark 3.6}. The above irreducible representation depends on
the three parameters $c\in \mbb{F}$ and $m_1,m_2\in\ol{1,n}$. It is
not highest-weight type because of the mixture of multiplication
operators and differential operators in (3.55), (3.56) and
(3.58)-(3.61). Since ${\msr B}$ is not completely reducible as a
module of the Lie subalgebra $\sum_{i,j=1}^m\mbb
F(E_{i,j}-E_{m+1+j,m+1+i})$ by \cite{LX1} when $m\geq 2$ and
$m_1<m$, $\msr A$  is not a unitary $sp(2m+2,\mbb{F})$-module.
Expression (3.62) shows that $\msr A$ is a weight
$sp(2m+2,\mbb{F})$-module with finite-dimensional weight
subspaces.\psp

{\it Case 4. $S\neq\emptyset,\;\vec a\neq 0, \;a_{i_0}\neq 0$ for
some $m+1+i_0\in S\bigcap\ol{m+2,2m+1}$ if
$S\bigcap\ol{m+2,2m+1}\neq\emptyset$, and $a_{m+1+j_0}\neq 0$ for
some $j_0\in S\bigcap\ol{1,m+1}$ if
$S\bigcap\ol{1,m+1}\neq\emptyset$.}\psp

We take (3.49)-(3.62). By the above assumption, $b_{j_0}\neq 0$ for
some $j_0\in\ol{1,m_1}$, and $a_{i_0}\neq 0$ for some
$i_0\in\ol{m_2+1,m}$ if $m_2<m$.  Recall the representation
$\pi_{c,S}^{\vec a}$ of $sp(2m+2,\mbb{F})$ is the representation
$\pi_{c,S}$ (cf. (3.52)-(3.62)) on $\msr A_{\vec a}$ (cf. (1.7)).
Under the assumption, we have the following fourth result in this
 section:\psp

{\bf Theorem 3.7}. {\it The representation $\pi_{c,S}^{\vec a}$ of
$sp(2m+2,\mbb{F})$ is irreducible for any $c\in\mbb{F}$.}

{\it Proof}.  Let ${\msr M}$ be a nonzero
$sp(2m+2,\mbb{F})$-submodule of ${\msr A}_{\vec a}$. Take any $0\neq
fe^{\vec a\cdot\vec x}\in \msr{M}$ with $f\in \msr{A}$.
 By the
assumption, $a_0=0$. Repeatedly applying the first equation in
(3.57) to $fe^{\vec a\cdot\vec x}$ if necessary, we may assume
$f\in\msr B=\mbb{F}[x_1,...,x_m,y_1,...,y_m]$. Then (3.59) yields
$$(E_{2m+2,i}+E_{m+1+i,m+1}+a_i)(fe^{\vec a\cdot\vec x})=-\ptl_{x_i}(f)e^{\vec a\cdot\vec
x}\in\msr M\;\;\for\;\;i\in\ol{m_1+1,m}.\eqno(3.77)$$ Moreover,
(3.58) yields
$$(E_{i,m+1}-E_{2m+2,m+1+i}-b_j)(fe^{\vec a\cdot\vec x})=\ptl_{y_j}(f)e^{\vec a\cdot\vec
x}\in\msr M\;\;\for\;\;j\in\ol{1,m_2}.\eqno(3.78)$$ Repeatedly
applying (3.77) and (3.78) if necessary, we can assume
$$f\in\mbb{F}[x_1,...,x_{m_1},y_{m_2+1},...,y_m].\eqno(3.79)$$

According to (3.55),
$$(E_{i,m+1+j_0}+E_{j_0,m+1+i}-a_{j_0}b_i-a_ib_{j_0})(fe^{\vec a\cdot\vec x})=(b_{j_0}\ptl_{x_i}+b_i\ptl_{x_{j_0}})(f)e^{\vec a\cdot\vec
x}\in\msr M\eqno(3.80)$$ for $i\in\ol{1,m_1}$. Taking $i=j_0$ in
(3.80), we get $\ptl_{x_{j_0}}(f)e^{\vec a\cdot\vec x}\in\msr M$.
Substituting it to (2.80) for general $i$, we obtain
$$\ptl_{x_i}(f)e^{\vec a\cdot\vec
x}\in\msr M\qquad\for\;\;i\in\ol{1,m_1}.\eqno(3.81)$$
 Moreover, (3.56) yields
$$(E_{m+1+j,i_0}+E_{m+1+i_0,j}-a_jb_{i_0}-a_{i_0}b_j)(fe^{\vec a\cdot\vec x})
=(a_{i_0}\ptl_{y_j}+a_j\ptl_{y_{i_0}})(f)e^{\vec a\cdot\vec
x}\in\msr M\eqno(3.82)$$ for $j\in\ol{m_2+1,m}.$ Letting $j=i_0$ in
(3.84), we find $\ptl_{y_{i_0}}(f)e^{\vec a\cdot\vec x}\in\msr M$.
Substituting it to (3.82) for general $j$, we get
$$\ptl_{y_j}(f)e^{\vec a\cdot\vec
x}\in\msr M\qquad\for\;\;j\in\ol{m_2+1,m}.\eqno(3.83)$$
 Repeatedly applying
(3.81) and (3.83) if necessary, we obtain $e^{\vec a\cdot\vec x}\in
\msr{M}$. Equivalently, $\msr{A}_{\vec a,0}\subset\msr{M}$ (cf.
(2.13)).

Suppose that for some $\ell\in\mbb{N}$, $\msr{A}_{\vec
a,k}\subset\msr{M}$ whenever $\ell\geq k\in\mbb{N}$. For any
$ge^{\vec a\cdot\vec x}\in \msr{A}_{\vec a,\ell}$, (3.59) implies
$$(E_{2m+2,i}+E_{m+1+i,m+1}-y_i\ptl_{x_0})(ge^{\vec a\cdot\vec
x})=x_ige^{\vec a\cdot\vec
x}\in\msr{M}\qquad\for\;\;i\in\ol{1,m_1}\eqno(3.84)$$ and (3.58)
leads to
$$(E_{2m+2,m+1+i}-E_{i,m+1}+x_i\ptl_{x_0})(ge^{\vec a\cdot\vec x})=y_jge^{\vec a\cdot\vec
x}\in\msr{M}\qquad\for\;\;j\in \ol{m_2+1,m}.\eqno(3.85)$$ Moreover,
(3.55) gives
$$(E_{i,m+1+j_0}+E_{j_0,m+1+i})(ge^{\vec a\cdot\vec
x})=[b_{j_0}x_i+x_i\ptl_{y_{j_0}}+(\ptl_{x_{j_0}}+a_{j_0})(\ptl_{y_i}+b_i)](g)e^{\vec
a\cdot\vec x}\in\msr{M}\eqno(3.86)$$ if $i\in\ol{m_1+1,m_2}$, and
$$(E_{i,m+1+j_0}+E_{j_0,m+1+i})(ge^{\vec a\cdot\vec
x})=[b_{j_0}x_i-a_{j_0}y_i+x_i\ptl_{y_{j_0}}-y_i\ptl_{x_{j_0}}](g)e^{\vec
a\cdot\vec x}\in\msr{M}\eqno(3.87)$$ if $i\in\ol{m_2+1,m}$. Note
that the inductional assumption imply
$$[x_i\ptl_{y_{j_0}}+(\ptl_{x_{j_0}}+a_{j_0})(\ptl_{y_i}+b_i)](g)e^{\vec
a\cdot\vec x}\in\msr{M}\eqno(3.88)$$  if $i\in\ol{m_1+1,m_2}$, and
$$[-a_{j_0}y_i+x_i\ptl_{y_{j_0}}-y_i\ptl_{x_{j_0}}](g)e^{\vec
a\cdot\vec x}\in\msr{M}\eqno(3.89)$$ by (3.87) if
$i\in\ol{m_2+1,m}$. Thus
$$x_ige^{\vec a\cdot\vec
x}\in\msr{M}\qquad\for\;\;i\in\ol{m_1+1,m}.\eqno(3.90)$$

On the other hand, (3.56) yields
$$(E_{m+1+j,i_0}+E_{m+1+i_0,j})(ge^{\vec a\cdot\vec
x})=(a_{i_0}y_j-b_{i_0}x_j+y_j\ptl_{x_{i_0}}-x_j\ptl_{y_{i_0}})(g)e^{\vec
a\cdot\vec x}\in\msr{M}\eqno(3.91)$$ if $j\in\ol{1,m_1}$, and
$$(E_{m+1+j,i_0}+E_{m+1+i_0,j})(ge^{\vec a\cdot\vec
x})=[a_{i_0}y_j+y_j\ptl_{x_{i_0}}+(\ptl_{x_j}+a_j)(\ptl_{y_{i_0}}+b_{i_0})](g)e^{\vec
a\cdot\vec x}\in\msr{M}\eqno(3.92)$$ if $j\in\ol{m_1+1,m_2}$.
Observe that the inductional assumption imply
$$(-b_{i_0}x_j+y_j\ptl_{x_{i_0}}-x_j\ptl_{y_{i_0}})(g)e^{\vec
a\cdot\vec x}\in\msr{M}\eqno(3.93)$$ by (3.84) if $j\in\ol{1,m_1}$,
and
$$[y_j\ptl_{x_{i_0}}+(\ptl_{x_j}+a_j)(\ptl_{y_{i_0}}+b_{i_0})](g)e^{\vec
a\cdot\vec x}\in\msr{M}\eqno(3.94)$$ if $j\in\ol{m_1+1,m_2}$. Hence
$$y_jge^{\vec a\cdot\vec x}\in\msr{M}\qquad\for\;\;j\in
\ol{1,m_2}.\eqno(3.95)$$  Moreover, (3.61) yields
\begin{eqnarray*}\qquad& &(E_{m+1,m+1+j_0}+E_{j_0,2m+2})(ge^{\vec a\cdot\vec
x})\\&=&[b_{j_0}x_0+x_0\ptl_{y_{j_0}}+(\td
D-\sum_{i=1}^{m_1}a_ix_i+\sum_{j=m_1+1}^ma_jx_j+\sum_{r=1}^{m_2}b_ry_r\\
& &-\sum_{s=m_2+1}^mb_sy_s +\td
c+1)(a_{j_0}+\ptl_{x_{j_0}})](g)e^{\vec a\cdot\vec
x}\in\msr{M}.\hspace{4.8cm}(3.96)\end{eqnarray*}  Note that
$$x_0\ptl_{y_{j_0}}(g)e^{\vec
a\cdot\vec x}\in\msr{A}_{\vec a,\ell}\subset\msr{M};\;\; (\td D +\td
c+1)(\ptl_{x_{j_0}}(g))e^{\vec a\cdot\vec x}\in\msr{A}_{\vec
a,\ell-1}\subset\msr{M}.\eqno(3.97)$$ Now (3.84), (3.85), (3.90) and
(3.95)-(3.97) imply $x_0ge^{\vec a\cdot\vec x}\in\msr{M}$.
Therefore, $\msr{A}_{\vec a,\ell+1}\subset\msr{M}$. By induction,
$\msr{A}_{\vec a,\ell}\subset\msr{M}$ for any $\ell\in\mbb{N}$. So
$\msr{A}_{\vec a}=\msr{M}$. Hence $\msr{A}_{\vec a}$ is an
irreducible $sp(2m+2,\mbb{F})$-module. $\qquad\Box$\psp

With respect to  the restricted representation $\pi_{c,S}^{\vec 0}$,
${\msr A}$ is an infinite-dimensional  weight
$sp(2m+2,\mbb{F})$-module with finite-dimensional weight subspaces
by (3.52)-(3.54) and (3.62). Now Theorem 2 follows from Theorems
3.2, 3.3. 3.5 and 3.7.

 \end{document}